\documentclass{amsart}
\usepackage{amssymb,amsmath,amscd,amsthm}
\usepackage{amsfonts,epsfig,latexsym,graphicx,amssymb,pict2e}

\newtheorem{Thm}{Theorem}
\newtheorem{Lm}{Lemma}[section]
\newtheorem{Prop}{Proposition}
\newtheorem{Def}[Lm]{Definition}
\newtheorem{Rem}[Lm]{Remark}
\newtheorem{Cor}{Corollary}

\def\bdef{\begin{Def}}
\def\endef{\end{Def}}
\def\bthm{\begin{Thm}}
\def\ethm{\end{Thm}}
\def\bprop{\begin{Prop}}
\def\enprop{\end{Prop}}
\def\blm{\begin{Lm}}
\def\elm{\end{Lm}}
\def\bcor{\begin{Cor}}
\def\ecor{\end{Cor}}
\def\brm{\begin{Rem}}
\def\erm{\end{Rem}}
\def\bfig{\begin{picture}}
\def\efig{\end{picture}}
\def\beq{\begin{eqnarray}}
\def\eneq{\end{eqnarray}}
\def\beal{\begin{aligned}}
\def\enal{\end{aligned}}

\title[Hyperbolicity of the Trace Map]{Hyperbolicity of the Trace Map for the Weakly Coupled Fibonacci Hamiltonian}

\author{David Damanik}

\address{Department of Mathematics, Rice University, Houston, TX~77005, USA}

\email{damanik@rice.edu}

\thanks{D.\ D.\ was supported in part by NSF grant
DMS--0653720.}

\author{Anton Gorodetski}

\address{Department of Mathematics, University of California, Irvine CA 92697, USA}

\email{asgor@math.uci.edu}

\begin{document}

\begin{abstract}
We consider the trace map associated with the Fibonacci
Hamiltonian as a diffeomorphism on the invariant surface
associated with a given coupling constant and prove that the
non-wandering set of this map is hyperbolic if the coupling is
sufficiently small. As a consequence, for these values of the
coupling constant, the local and global Hausdorff dimension and
the local and global box counting dimension of the spectrum of the
Fibonacci Hamiltonian all coincide and are smooth functions of the
coupling constant.
\end{abstract}

\maketitle


\section{Introduction}

The Fibonacci Hamiltonian is the most prominent model in the study
of electronic properties of quasicrystals. It is given by the
discrete one-dimensional Schr\"odinger operator
$$
[H_{V,\omega} u](n) = u(n+1) + u(n-1) + V \chi_{[1-\alpha,1)}(n
\alpha + \omega \!\!\!\! \mod 1) u(n),
$$
where $V > 0$ is the coupling constant, $\alpha =
\frac{\sqrt{5}-1}{2}$ is the frequency, and $\omega \in [0,1)$ is
the phase.

This operator family displays a number of interesting phenomena,
such as Cantor spectrum of zero Lebesgue measure \cite{S89} and
purely singular continuous spectral measure for all phases
\cite{DL}. Moreover, it was recently shown that is also gives rise
to anomalous transport \cite{DT}. We refer the reader to the
survey articles \cite{D00,D07,S95} for further information and
references.

Already the earliest papers on this model, \cite{KKT,OPRSS},
realized the importance of a certain renormalization procedure in
its study. This led in particular to a consideration of the
following dynamical system, the so-called trace map,
$$
T:\Bbb{R}^3\to \Bbb{R}^3, \; T(x,y,z)=(2xy-z,x,y),
$$
whose properties are closely related to all the spectral
properties mentioned above. The existence of the trace map and its
connection to spectral properties of the operators is a
consequence of the invariance of the potential under a
substitution rule. This works in great generality; see the surveys
mentioned above and references therein. In the Fibonacci
case, the existence of a first integral is an additional useful property,
which allows one to restrict $T$ to invariant surfaces.

Fix some coupling constant $V$. For a complete spectral study of
the operator family $\{ H_{V,\omega} \}_{\omega \in [0,1)}$, it
suffices to study $T$ on a single invariant surface $S_V$. This
phenomenon is a peculiarity of the choice of the model (a discrete
Schr\"odinger operator in math terminology or an on-site lattice
model in physics terminology) and does not follow solely from the
symmetries coming from the invariance under the Fibonacci
substitution. For example, in continuum analogs of the Fibonacci
Hamiltonian, the invariant surface will in general be energy-dependent.

Let us denote the restriction of $T$ to the invariant surface
$S_V$ by $T_V$. As we will discuss in more detail below, it is of
interest to study the non-wandering set of this surface
diffeomorphism because it is closely connected to the spectrum of
$H_{V,\omega}$.\footnote{By a strong convergence argument, it
follows that the spectrum of $H_{V,\omega}$ does not depend on
$\omega$. It does, however, depend on $V$.} This correspondence in
turn allows one to show that the spectrum has zero Lebesgue
measure. It is then natural to investigate its fractal dimension.
A number of papers have studied this problem; for example,
\cite{DEGT,LW,Ra}. As pointed out in \cite{DEGT}, the work of
Casdagli, \cite{Cas}, has very important consequences for the
fractal dimension of the spectrum as a function of $V$. Casdagli
studied the map $T_V$ and proved, for $V \ge 16$, that the
non-wandering set is hyperbolic. Combining this with results in
hyperbolic dynamics, it follows that the local and global
Hasudorff and box counting dimensions of the spectrum all coincide
and are smooth functions of $V$. This result was crucial for the
work \cite{DEGT}, which determined the exact asymptotic behavior
of this function of $V$ as $V$ tends to infinity. It was shown
that
$$
\lim_{V \to \infty} \dim \sigma(H_{V,\omega}) \cdot \log V = \log
(1 + \sqrt{2}).
$$

Of course, the asymptotic behavior of the dimension of the
spectrum as $V$ approaches zero is of interest as well. Given the
discussion above, the natural first step is to prove the analogue
of Casdagli's result at small coupling. This is exactly what we do
in this paper. We will show that, for $V$ sufficiently small, the
non-wandering set of $T_V$ is hyperbolic and hence we obtain the
same consequences for the dimension of the spectrum as those
mentioned above in this coupling regime.

The structure of the paper is as follows.
Section~\ref{s.background} gives a more explicit description of
the previous results on the Fibonacci trace map, recalls some
useful general results from hyperbolic dynamics, and states the
main result of the paper --- the hyperbolicity of the
non-wandering set of the trace map for sufficiently small coupling
$V$. Sections~\ref{s.vequzero}--\ref{s.ppe} contain the proof of
this result. More precisely, Section~\ref{s.vequzero} contains a
discussion of the case $V=0$, Section~\ref{s.s04} studies the
dynamics of the trace map near a singular point and formulates the
crucial Proposition~\ref{p.expanding}. Section~\ref{s.s05} shows
how the main result follows from it, and finally,
Section~\ref{s.ppe} contains a proof of
Proposition~\ref{p.expanding}.

After this paper was finished we learned that Serge Cantat
provided a proof of uniform hyperbolicity of the trace map for all
non-zero values of the coupling constant \cite{Can}. Our results
were obtained independently and we use completely different
methods.

\bigskip

\noindent\textit{Acknowledgments.} We wish to express our
gratitude to Helge Kr\"uger, John Roberts, and Amie Wilkinson for
useful discussions.

\section{Background and Main Result}\label{s.background}

In this section we expand on the introduction and state
definitions and previous results more carefully. This will
eventually lead us to the statement of our main result in
Theorem~\ref{t.main} below.

\subsection{Description of the Trace Map and Previous Results}

The main tool that we are using here is the so called
\textit{trace map}. It was originally introduced in \cite{K,KKT};
further useful references include \cite{BGJ,BR,HM,Ro}. Let us quickly
recall how it arises from the substitution invariance of the
Fibonacci potential; see \cite{S87} for detailed proofs of some of
the statements below.

The one step transfer matrices associated with the difference
equation $H_{V,\omega} u = E u$ are given by
$$
T_{V,\omega}(m,E) = \begin{pmatrix} E - V \chi_{[1-\alpha,1)}(m
\alpha + \omega \!\!\!\! \mod 1) & -1 \\
1 & 0 \end{pmatrix}.
$$
Denote the Fibonacci numbers by $\{F_k\}$, that is, $F_0 = F_1 =
1$ and $F_{k+1} = F_k + F_{k_1}$ for $k \ge 1$. Then, one can show
that the matrices
$$
M_{-1}(E) = \begin{pmatrix} 1 & -V \\ 0 & 1
\end{pmatrix} , \quad M_0(E) =
\begin{pmatrix} E & -1 \\ 1 & 0 \end{pmatrix}
$$
and
$$
M_k(E) = T_{V,0}(F_k,E) \times \cdots \times T_{V,0}(1,E) \quad
\text{ for } k \ge 1
$$
obey the recursive relations
$$
M_{k+1}(E) = M_{k-1}(E) M_k(E)
$$
for $k \ge 0$. Passing to the variables
$$
x_k(E) = \frac12 \mathrm{Tr} M_k(E),
$$
this in turn implies
$$
x_{k+1}(E) = 2 x_k(E) x_{k-1}(E) - x_{k-2}(E).
$$
These recursion relations exhibit a conserved quantity; namely, we
have
$$
x_{k+1}(E)^2+x_k(E)^2+x_{k-1}(E)^2 - 2 x_{k+1}(E) x_k(E)
x_{k-1}(E) -1 = \frac{V^2}{4}
$$
for every $k \ge 0$.

\begin{figure}\begin{minipage}[t]{5cm} \includegraphics[width=1.2\textwidth]{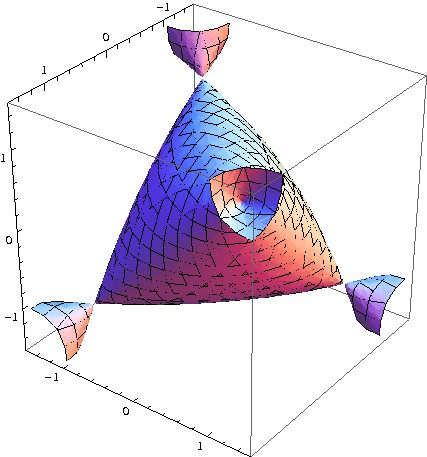} \par
\caption{The surface $S_{0.01}$.} \label{fig:s0.01}
\end{minipage} \hfill \begin{minipage}[t]{5cm} \includegraphics[width=1.2\textwidth]{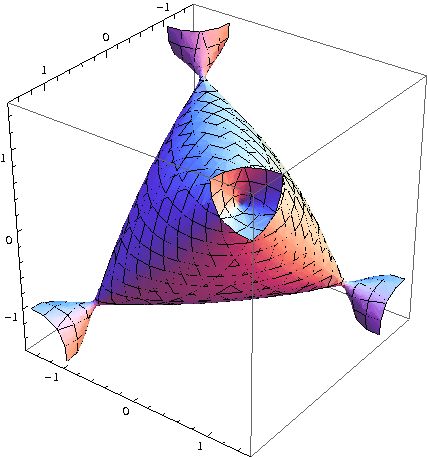} \par
\caption{The surface $S_{0.1}$.} \label{fig:s0.1}
\end{minipage} \end{figure}

\begin{figure} \begin{minipage}[t]{5cm} \includegraphics[width=1.2\textwidth]{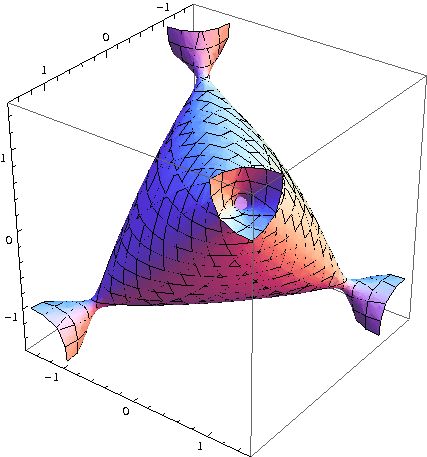} \par
\caption{The surface $S_{0.2}$.} \label{fig:s0.2}
\end{minipage} \hfill \begin{minipage}[t]{5cm} \includegraphics[width=1.2\textwidth]{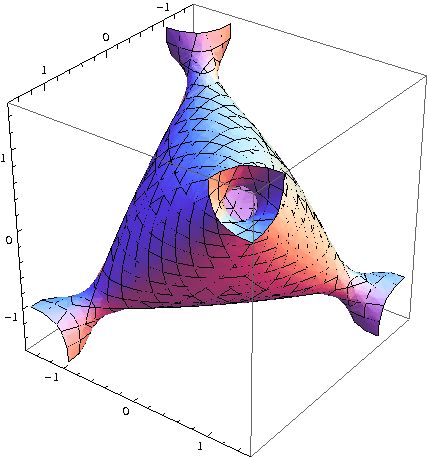}
\par \caption{The surface $S_{0.5}$.} \label{fig:s0.5} \end{minipage}
\end{figure}

Given these observations, it is then convenient to introduce the
\textit{trace map}
$$
T : \Bbb{R}^3 \to \Bbb{R}^3, \; T(x,y,z)=(2xy-z,x,y).
$$
The following function
$$
G(x,y,z) = x^2+y^2+z^2-2xyz-1
$$
is invariant under the action of $T$,\footnote{The function
$G(x,y,z)$ is called the \textit{Fricke character}.} and hence $T$
preserves the family of cubic surfaces\footnote{The surface $S_0$
is called the \textit{Cayley cubic}.}
$$
S_V = \left\{(x,y,z)\in \Bbb{R}^3 : x^2+y^2+z^2-2xyz=1+
\frac{V^2}{4} \right\}.
$$
Plots of the surfaces $S_{0.01}$, $S_{0.1}$, $S_{0.2}$, and
$S_{0.5}$ are given in Figures~\ref{fig:s0.01}--\ref{fig:s0.5},
respectively.

Denote by $\ell_V$ the line
$$
\ell_V = \left\{ \left(\frac{E-V}{2}, \frac{E}{2}, 1 \right) : E
\in \Bbb{R} \right\}.
$$
It is easy to check that $\ell_V \subset S_V$.

S\"ut\H{o} proved the following central result in \cite{S87}.

\bthm[S\"ut\H{o} 1987]\label{spectrum}  An energy $E$ belongs to
the spectrum of $H_{V,\omega}$ if and only if the positive
semiorbit of the point $(\frac{E-V}{2}, \frac{E}{2}, 1)$ under
iterates of the trace map $T$ is bounded. \ethm

It is of course natural to consider the restriction $T_{V}$ of the
trace map $T$ to the invariant surface $S_V$. That is, $T_{V}:S_V
\to S_V$, $T_{V}=T|_{S_V}$. Denote by $\Omega_{V}$ the set of
points in $S_V$ whose full orbits under $T_{V}$ are bounded. {\it
A priori} the set of bounded orbits of $T_V$ could be different
from the non-wandering set\footnote{A point $p\in M$ of a
diffeomorphism $f:M\to M$ is wandering if there exists a
neighborhood $O(p)\subset M$ such that $f^k(O)\cap O = \emptyset$
for any $k\in \mathbb{Z}\backslash 0$. The non-wandering set of
$f$ is the set of points that are not wandering.} of $T_V$, but
our construction of the Markov partition and analysis of the
behavior of $T_V$ near singularities show that in our case these
two sets do coincide. Notice that this is parallel to the
construction of the symbolic coding in \cite{Cas}.

Let us recall that an invariant closed set $\Lambda$ of a
diffeomorphism $f : M \to M$ is \textit{hyperbolic} if there
exists a splitting of a tangent space $T_xM=E^u_x\oplus E^u_x$ at
every point $x\in \Lambda$ such that this splitting is invariant
under $Df$ and the differential $Df$ exponentially contracts
vectors from stable subspaces $\{E^s_x\}$ and exponentially
expands vectors from unstable subspaces $\{E^u_x\}$. A hyperbolic
set $\Lambda$ of a diffeomorphism $f : M \to M$ is \textit{locally
maximal} if there exists a neighborhood $U(\Lambda)$ such that
$$
\Lambda=\bigcap_{n\in\Bbb{Z}}f^n(U).
$$

The second central result about the trace map we wish to recall is
due to Casdagli; see \cite{Cas}.

\bthm[Casdagli 1986]\label{Casdagli} For $V\ge 16$, the set
$\Omega_{V}$ is a locally maximal hyperbolic set of
$T_{V}:S_V \to S_V$. It is homeomorphic to a Cantor set. \ethm

\subsection{Some Properties of Locally Maximal Hyperbolic Invariant
Sets of Surface Diffeomorphisms}

Given Theorem~\ref{Casdagli}, several general results apply to the
trace map of the strongly coupled Fibonacci Hamiltonian. Let us
recall some of these results that yield interesting spectral
consequences, which are discussed below.

Consider a locally maximal invariant transitive hyperbolic set
$\Lambda\subset M$, $\dim M=2$, of a diffeomorphism $f\in \text{\,
\rm Diff\,}^r(M), r\ge 1$. We have
$\Lambda=\bigcap_{n\in\Bbb{Z}}f^n(U(\Lambda))$ for some neighborhood
$U(\Lambda)$. Assume also that $\dim E^u=\dim E^s=1$. Then, the
following properties hold.

\subsubsection{Stability}

There is a neighborhood $\mathcal{U}\subset \text{\, \rm
Diff\,}^1(M)$ of the map $f$ such that for every $g\in
\mathcal{U}$, the set $\Lambda_g = \bigcap_{n\in
\Bbb{Z}}g(U(\Lambda))$ is a locally maximal invariant hyperbolic
set set of $g$. Moreover, there is a homeomorphism $h:\Lambda\to
\Lambda_g$ that conjugates $f|_{\Lambda}$ and $g|_{\Lambda_g}$,
that is, the following diagram commutes:
$$
\CD \Lambda @>{f|_{\Lambda}}>>\Lambda \\
@V{h}VV@VV{h}V\\
\Lambda_g @>g|_{\Lambda_g}>>\Lambda_g\endCD
$$

\subsubsection{Invariant Manifolds}

For $x\in \Lambda$ and small $\varepsilon>0$, consider the local
stable and unstable sets
$$
W^s_{\varepsilon}(x)=\{w\in M : d(f^n(x), f^n(w))\le \varepsilon \
\ \text{\rm for all}\ \  n\ge 0\},
$$
$$
W^u_{\varepsilon}(x)=\{w\in M : d(f^n(x), f^n(w))\le \varepsilon \
\ \text{\rm for all}\ \  n\le 0\}.
$$
If $\varepsilon>0$ is small enough, these sets are embedded
$C^r$-disks with $T_xW^s_{\varepsilon}(x)=E^s_x$ and
$T_xW^u_{\varepsilon}(x)=E^u_x$. Define the (global) stable and
unstable sets as
$$
W^s(x)=\bigcup_{n\in \Bbb{N}}f^{-n}(W^s_{\varepsilon}(x)), \ \ \
W^u(x)=\bigcup_{n\in \Bbb{N}}f^{n}(W^u_{\varepsilon}(x)).
$$
Define also
$$
W^s(\Lambda)=\bigcup_{x\in \Lambda}W^s(x)\ \ \ \text{\rm and}\ \ \
W^u(\Lambda)=\bigcup_{x\in \Lambda}W^u(x).
$$

\subsubsection{Invariant Foliations}

A stable foliation for $\Lambda$ is a foliation $\mathcal{F}^s$ of a neighborhood of $\Lambda$ such that
\begin{itemize}
\item[(a)] for each $x\in \Lambda$, $\mathcal{F}(x)$, the leaf
containing $x$, is tangent to $E^s_x$,

\item[(b)] for each $x$ sufficiently close to $\Lambda$,
$f(\mathcal{F}^s(x))\subset \mathcal{F}^s(f(x))$.
\end{itemize}
An unstable foliation $\mathcal{F}^u$ can be defined in a similar way.

For a locally maximal hyperbolic set $\Lambda\subset M$ of a
$C^1$-diffeomorphism $f:M\to M$, $\dim M=2$, stable and unstable
$C^0$ foliations with $C^1$-leaves can be constructed \cite{M}. In
the case of $C^2$-diffeomorphism,  $C^1$ invariant foliations
exist (see, for example, \cite{PT}, Theorem 8 in Appendix 1).

\subsubsection{Local Hausdorff Dimension and Box Counting
Dimension}

Consider, for $x\in \Lambda$ and small $\varepsilon>0$, the set
$W^u_{\varepsilon}(x)\cap \Lambda$. Its Hausdorff dimension does
not depend on $x\in \Lambda$ and $\varepsilon>0$, and coincides
with its box counting dimension (see \cite{MM,T}):
$$
\dim_H W^u_{\varepsilon}(x)\cap \Lambda = \dim_B
W^u_{\varepsilon}(x)\cap \Lambda.
$$
In a similar way,
$$
\dim_H W^s_{\varepsilon}(x)\cap \Lambda= \dim_B
W^s_{\varepsilon}(x)\cap \Lambda.
$$
Denote $h^s = \dim_H W^s_{\varepsilon}(x)\cap \Lambda$ and $h^u =
\dim_H W^u_{\varepsilon}(x)\cap \Lambda$. We will say that $h^s$
and $h^u$ are the \textit{local stable and unstable Hausdorff
dimensions} of $\Lambda$.

For properly chosen small $\varepsilon > 0$, the sets
$W^u_{\varepsilon}(x)\cap \Lambda$ and $W^s_{\varepsilon}(x)\cap
\Lambda$ are dynamically defined Cantor sets (see \cite{PT1} for
definitions and proof), and this implies, in particular, that
$$
h^s<1 \ \ \  \text{\rm and} \ \ \ h^u<1,
$$
see, for example, Theorem 14.5 in \cite{P}.

\subsubsection{Global Hausdorff Dimension}

The Hausdorff dimension of $\Lambda$ is equal to its box
counting dimension, and
$$
\dim_H \Lambda=\dim_B \Lambda = h^s + h^u;
$$
see \cite{MM,PV}.

\subsubsection{Continuity of the Hausdorff Dimension}

The local Hausdorff dimensions $h^s(\Lambda)$ and $h^u(\Lambda)$
depend continuously on $f:M\to M$ in the $C^1$-topology; see
\cite{MM,PV}. Therefore, $\dim_H \Lambda_f = \dim_B \Lambda_f =
h^s(\Lambda_f) + h^u(\Lambda_f)$ also depends continuously on $f$
in the $C^1$-topology. Moreover, for a $C^r$ diffeomorphism
$f:M\to M$, $r\ge2$, the Hausdorff dimension of a hyperbolic set
$\Lambda_f$ is a $C^{r-1}$ function of $f$; see \cite{Ma}.

\brm For hyperbolic sets in dimension greater than two, many of
these properties do not hold in general; see \cite{P} for more
details. \erm

\subsection{Implications for the Trace Map and the Spectrum}

Due to Theorem~\ref{Casdagli}, for every $V\ge 16$, all the
properties from the previous subsection can be applied to the
hyperbolic set $\Omega_{\lambda}$ of the trace map $T_{V}:S_V \to
S_V$.

Moreover, the results in \cite[Section 2]{Cas} imply the following
statement. \blm\label{l.translargecoupling} For $V\ge 16$ and
every $x\in \Omega_{V}$, the stable manifold $W^s(x)$ intersects
the line $\ell_V$ transversally. \elm

The existence of a $C^1$-foliation $\mathcal{F}^s$ allows one to
locally consider the set $W^s(\Omega_{V})\cap \ell_V$ as a
$C^1$-image of the set $W^u_{\varepsilon}(x)\cap\Omega_{V}$. Due
to Theorem \ref{spectrum}, this implies the following properties
of the spectrum $\sigma(H_{V,\omega})$ for $V\ge 16$:

\begin{Cor}\label{c.specprop}
For $V\ge 16$, the following statements hold:

\begin{itemize}

\item[{\rm (i)}] The spectrum $\sigma(H_{V,\omega})$ depends
continuously on $V$ in the Hausdorff metric.

\item[{\rm (ii)}] For every small $\varepsilon>0$ and every $x\in
\sigma(H_{V,\omega})$, we have
\begin{align*}
\dim_H \left( (x-\varepsilon, x+\varepsilon) \cap
\sigma(H_{V,\omega}) \right) & = \dim_B \left( (x-\varepsilon,
x+\varepsilon) \cap \sigma(H_{V,\omega}) \right) \\ & = \dim_H
\sigma(H_{V,\omega}) \\ & =\dim_B \sigma(H_{V,\omega}).
\end{align*}

\item[{\rm (iii)}] The Hausdorff dimension $\dim_H
\sigma(H_{V,\omega})$ is a $C^{\infty}$-function of $V$, and is
strictly smaller than one.

\end{itemize}
\end{Cor}

\subsection{Hyperbolicity of the Trace Map for
Small Coupling}

We are now in a position to state the main result of this paper.

\bthm\label{t.main} There exists $V_0>0$ such that for every $V\in
(0,V_0)$, the following properties hold.

\begin{itemize}

\item[{\rm (i)}] The non-wandering set $\Omega_{V}\subset S_{V}$
of the map $T_{V}:S_{V}\to S_{V}$ is hyperbolic.

\item[{\rm (ii)}] The non-wandering set $\Omega_{V}\subset S_{V}$
is homeomorphic to a Cantor set, and $T_{V}|_{\Omega_{V}}$ is
conjugated to a topological Markov chain $\sigma_{A}:\Sigma^6_A\to
\Sigma^6_A$ with the matrix
$$
A=\begin{pmatrix}
    0 & 0 & 0 & 1 & 1 & 1 \\
    0 & 0 & 1 & 0 & 1 & 1 \\
    0 & 0 & 0 & 0 & 1 & 0 \\
    0 & 0 & 0 & 0 & 0 & 1 \\
    1 & 0 & 0 & 0 & 0 & 0 \\
    0 & 1 & 0 & 0 & 0 & 0 \\
  \end{pmatrix}.
$$

\item[{\rm (iii)}] For every $x\in \Omega_{V}$, the stable
manifold $W^s(x)$ intersects the line $\ell_V$ transversally.

\end{itemize}
\ethm

As before, we obtain the following consequences:

\begin{Cor}\label{c.specpropsl}
With $V_0 > 0$ from Theorem~\ref{t.main}, the following statements
hold for $V \in (0,V_0)$:

\begin{itemize}

\item[{\rm (i)}] The spectrum $\sigma(H_{V,\omega})$ depends
continuously on $V$ in the Hausdorff metric.

\item[{\rm (ii)}] For every small $\varepsilon>0$ and every $x\in
\sigma(H_{V,\omega})$, we have
\begin{align*}
\dim_H \left( (x-\varepsilon, x+\varepsilon) \cap
\sigma(H_{V,\omega}) \right) & = \dim_B \left( (x-\varepsilon,
x+\varepsilon) \cap \sigma(H_{V,\omega}) \right) \\ & = \dim_H
\sigma(H_{V,\omega}) \\ & =\dim_B \sigma(H_{V,\omega}).
\end{align*}

\item[{\rm (iii)}] The Hausdorff dimension $\dim_H
\sigma(H_{V,\omega})$ is a $C^{\infty}$-function of $V$, and is strictly smaller than one.

\end{itemize}
\end{Cor}

We expect these properties to be of similar importance in a study
of the asymptotic behavior of the fractal dimension of the
spectrum as $V \to 0$ as was the case in the large coupling
regime.

\section{Properties of the Trace Map for $V=0$}\label{s.vequzero}

Up to this point, we only considered the case $V > 0$. Since we
will regard the case of small positive $V$ as a small perturbation
of the case $V = 0$, we will also include the latter case in our
considerations. In fact, this section is devoted to the study of
this ``unperturbed case.''

Denote by $\mathbb{S}$ the part of the surface $S_0$ inside of the
cube $\{ |x|\le 1, |y|\le 1, |z|\le 1\}$. The surface $\mathbb{S}$
is homeomorphic to $S^2$, invariant, smooth everywhere except at the four
points $P_1=(1,1,1)$, $P_1=(1,-1,-1)$, $P_1=(-1,1,-1)$, and
$P_1=(-1,-1,1)$, where $\mathbb{S}$ has conic singularities,  and
the trace map $T$ restricted to $\mathbb{S}$ is a factor of a
hyperbolic automorphism of a two torus:
$$
\mathcal{A}(\theta, \varphi)=(\theta + \varphi, \theta)\ (\text{\rm
mod}\ 1).
$$
The semiconjugacy is given by the map
$$
F: (\theta, \varphi) \mapsto (\cos 2\pi(\theta + \varphi), \cos 2\pi
\theta, \cos 2\pi \varphi).
$$
The map $\mathcal{A}$ is hyperbolic, and is given by a matrix $A=\begin{pmatrix}
                                                                   1 & 1 \\
                                                                   1 & 0 \\
                                                                 \end{pmatrix}$
with eigenvalues
$$
\mu=\frac{1+\sqrt{5}}{2}\ \ \text{\rm and} \ \ \
-\mu^{-1}=\frac{1-\sqrt{5}}{2}.
$$
Let us denote by $\mathbf{v}^u, \mathbf{v}^u \in \mathbb{R}^2 $ the
unstable and stable eigenvectors of $A$:
$$
A\mathbf{v}^u=\mu\mathbf{v}^u, \
A\mathbf{v}^s=-\mu^{-1}\mathbf{v}^s, \
\|\mathbf{v}^u\|=\|\mathbf{v}^s\|=1.
$$
Fix some small $\zeta > 0$ and define the stable (resp., unstable)
cone fields on $\mathbb{R}^2$ in the following way:
\begin{align}\label{e.cones}
    K^s_p & =\{\mathbf{v}\in T_p\mathbb{R}^2 :
\mathbf{v}=v^u\mathbf{v}^u+v^s\mathbf{v}^s, \ |v^s| > \zeta^{-1}
|v^u|\},\\
\nonumber    K^u_p & =\{\mathbf{v}\in T_p\mathbb{R}^2 :
\mathbf{v}=v^u\mathbf{v}^u+v^s\mathbf{v}^s, \ |v^u|
> \zeta^{-1} |v^s|\}.
\end{align}
These cone fields are invariant:
$$
\forall\ \mathbf{v}\in K^u_p\ \ \ \ A\mathbf{v}\in K^u_{A(p)},
$$
$$
\forall\ \mathbf{v}\in K^s_p\ \ \ \ A^{-1}\mathbf{v}\in
K^s_{A^{-1}(p)}.
$$
Also, the iterates of the map $A$ expand vectors from the unstable
cones, and the iterates of the map $A^{-1}$ expand vectors from the
stable cones:
$$
\forall\ \mathbf{v}\in K^u_p\ \ \ \ \forall\ n\in \mathbb{N}\ \ \
\ |A^n\mathbf{v}| > \frac{1}{\sqrt{1 + \zeta^2}}
\lambda^n|\mathbf{v}|,
$$
$$
\forall\ \mathbf{v}\in K^s_p\ \ \ \ \forall\ n\in \mathbb{N}\ \ \
\ |A^{-n}\mathbf{v}| > \frac{1}{\sqrt{1 + \zeta^2}}
\lambda^n|\mathbf{v}|.
$$
The families of cones $\{K^s\}$ and $\{K^u\}$ invariant under
$\mathcal{A}$ can be also considered on $\mathbb{T}^2$.

The differential of the semiconjugacy $F$ sends these cone families
to stable and unstable cone families on $\mathbb{S}\backslash \{P_1,
P_2, P_3, P_4\}$. Let us denote these images by $\{\mathcal{K}^s\}$
and $\{\mathcal{K}^u\}$.

\blm\label{l.differential} The differential of the semiconjugacy $DF$ induces a map of
the unit bundle of $\mathbb{T}^2$ to the unit bundle of
$\mathbb{S}\backslash \{P_1, P_2, P_3, P_4\}$. The derivatives of
the restrictions of this map to a fiber are uniformly bounded. In
particular, the sizes of cones in families $\{\mathcal{K}^s\}$ and
$\{\mathcal{K}^u\}$ are uniformly bounded away from zero.  \elm

\begin{proof}
Choose small enough neighborhoods $U_1(P_1), U_2(P_2), U_3(P_3)$,
and $U_4(P_4)$ in $\mathbb{S}$. The complement
$$
\hat{\mathbb{S}} = \mathbb{S} \setminus \Big( \bigcup_{j=1}^4
U_j(P_j) \Big)
$$
is compact, so $F^{-1}(\hat{\mathbb{S}})$ is also compact, and the
action of $DF$ on a fiber of the unit bundle over points of
$F^{-1}(\hat{\mathbb{S}})$ has uniformly bounded derivatives.

Due to the symmetries of the trace map (see, for example, \cite{K}
for the detailed description of the symmetries of the trace map)
and the semiconjugacy $F$, it is enough to consider a neighborhood
$U_1(P_1)$. Hence, it is enough to consider the map $F$ in a
neighborhood of the point $(0,0)$.

The differential of $F$ has the form
$$
DF=-2\pi\begin{pmatrix}
          \sin 2\pi(\theta+\varphi) & \sin 2\pi(\theta+\varphi) \\
          \sin 2\pi\theta & 0 \\
          0 & \sin 2\pi\varphi \\
        \end{pmatrix}.
$$
If $(\theta, \varphi)$ is in a small neighborhood of $(0,0)$, then
$$
DF\sim -4\pi^2\begin{pmatrix}
          \theta+\varphi & \theta+\varphi \\
        \theta & 0 \\
          0 & \varphi \\
        \end{pmatrix}.
$$

Therefore, up to a multiplicative constant and higher order terms,
the image of the vector $(1,0)$ under $DF$ is $(\theta+\varphi,
\theta, 0)$, and the image of the vector $(0,1)$ is
$(\theta+\varphi, 0, \varphi)$. In order to estimate the
derivative of the projective action of $DF$ on a fiber over a
point $(\theta, \varphi)$, it is enough to estimate the angle
between images of basis vectors, and the ratio of the lengths of
the images of these vectors.

If $\alpha$ is the angle between $DF_{(\theta, \varphi)}(1,0)$ and
$DF_{(\theta, \varphi)}(0,1)$, then
$$
\cos \alpha \sim
\frac{(\theta+\varphi)^2}{\sqrt{(\theta+\varphi)^2+\theta^2}\cdot
\sqrt{(\theta+\varphi)^2+\varphi^2}}=\frac{1}{\sqrt{1+\frac{\theta^2+\varphi^2}{(\theta+\varphi)^2}
+\frac{\theta^2\varphi^2}{(\theta+\varphi)^4}}}.
$$
We have
$$
\frac{\theta^2+\varphi^2}{(\theta+\varphi)^2}\ge \frac{1}{2},\ \ \
\text{\rm and}\ \ \
\frac{\theta^2\varphi^2}{(\theta+\varphi)^4}\ge 0,
$$
so $\cos \alpha \le \sqrt{\frac{2}{3}}+0.001<1$ if $(\theta,
\varphi)$ is close enough to $(0,0)$.

Now let us estimate the ratio of the lengths of $DF_{(\theta,
\varphi)}(1,0)$ and $DF_{(\theta, \varphi)}(0,1)$. Up to higher
order terms it is equal to
$$
\frac{\sqrt{\varphi^2+2\theta^2+2\theta\varphi}}{\sqrt{2\varphi^2+\theta^2+2\theta\varphi}}=\sqrt{\frac{1+2t^2+2t}{2+t^2+2t}}=
\sqrt{\frac{2(t+\frac{1}{2})^2+\frac{1}{2}}{(t+1)^2+1}},
$$
where $t=\frac{\theta}{\varphi}\in \mathbb{R} \cup{\{\infty\}}$,
and this function is bounded from above and bounded away from
zero. Lemma \ref{l.differential} is proved.
\end{proof}

\section{The Structure of the Trace Map in a Neighborhood of a Singular
Point}\label{s.s04}

Due to the symmetries of the trace map it is enough to consider
the dynamics of $T$ in a neighborhood of $P_1=(1,1,1)$. Let
$U\subset \mathbb{R}^3$ be a small neighborhood of $P_1$ in
$\mathbb{R}^3$. Let us consider the set $Per_2(T)$ of periodic
points of $T$ of period 2.

\begin{figure} \begin{minipage}[t]{5cm} \includegraphics[width=1.2\textwidth]{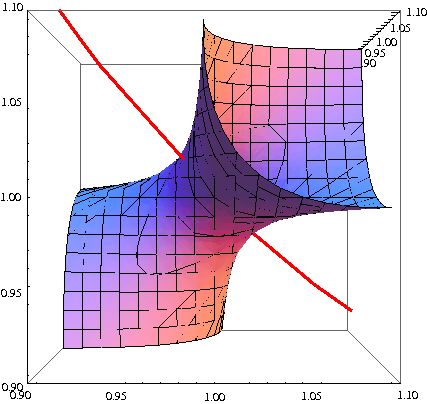} \par
\caption{$S_{0.1}$ and $Per_2(T)$ near $(1,1,1)$.}
\label{fig:sing0.0} \end{minipage} \hfill
\begin{minipage}[t]{5cm}
\includegraphics[width=1.2\textwidth]{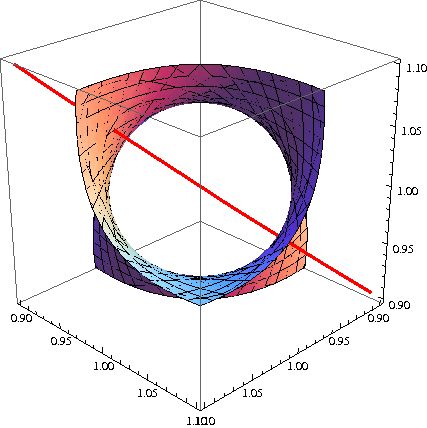}
\par \caption{$S_{0.2}$ and $Per_2(T)$ near $(1,1,1)$.} \label{fig:sing0.2}
\end{minipage} \end{figure}

\blm\label{l.periodtwo} We have
$$
Per_2(T)=\left\{(x,y,z) : x \in (-\infty, \tfrac{1}{2})\cup
(\tfrac{1}{2}, \infty), \ y=\frac{x}{2x-1}, \ z=x\right\}.
$$
\elm
\begin{proof}
Direct calculation.
\end{proof}
Notice that in a neighborhood $U$ the intersection $I\equiv
Per_2(T)\cap U$ is a smooth curve that is a normally hyperbolic
with respect to $T$ (see, for example, \cite{PT}, Appendix 1, for
the formal definition of normal hyperbolicity). Therefore, the
local center-stable manifold  $W_{loc}^{cs}(I)$ and the local
center-unstable manifold $W^{cu}_{loc}(I)$ defined by
$$
W_{loc}^{cs}(I)=\left\{p\in U : T^n(p) \in U \ \text{\rm for all}\
n\in \mathbb{N}\right\},
$$
$$
W_{loc}^{cu}(I)=\left\{p\in U : T^{-n}(p) \in U \ \text{\rm for
all}\ n\in \mathbb{N}\right\}
$$
are smooth two-dimensional surfaces. Also, the local strong stable
manifold $W^{ss}_{loc}(P_1)$ and the local strong unstable
manifold $W^{uu}_{loc}(P_1)$ of the fixed point $P_1$, defined by
$$
W^{ss}_{loc}(P_1)=\left\{p\in W^{cs}_{loc}(I) : T^n(p)\to P_1 \
\text{\rm as}\ n\to +\infty\right\},
$$
$$
W^{uu}_{loc}(P_1)=\left\{p\in W^{cu}_{loc}(I) : T^{-n}(p)\to P_1 \
\text{\rm as}\ n\to +\infty\right\},
$$
are smooth curves.

Let $\Phi:U\to \mathbb{R}^3$ be a smooth change of coordinates such that $\Phi(P_1)=(0,0,0)$ and
\begin{itemize}
\item
$\Phi(I)$ is a part of the line $\{x=0, z=0\}$;
\item
$\Phi(W^{cs}_{loc}(I))$ is a part of the plane $\{z=0\}$;
\item
$\Phi(W^{cu}_{loc}(I))$ is a part of the plane $\{x=0\}$;
\item
$\Phi(W^{ss}_{loc}(P_1))$ is a part of the line $\{y=0, z=0\}$;
\item
$\Phi(W^{uu}_{loc}(P_1))$ is a part of the line $\{x=0, y=0\}$.
\end{itemize}
Denote $f=\Phi\circ T\circ \Phi^{-1}$.

In this case,
$$
Df(0,0,0)=D(\Phi\circ T\circ \Phi^{-1})(0,0,0)=\begin{pmatrix}
                                       \lambda^{-1} & 0 & 0 \\
                                       0 & -1 & 0 \\
                                       0 & 0 & \lambda \\
                                     \end{pmatrix},
$$
where $\lambda$ is a largest eigenvalue of the differential
$DT(P_1):T_{P_1}\mathbb{R}^3\to T_{P_1}\mathbb{R}^3$,
$$
DT(P_1)=\begin{pmatrix}
                                       2 & 2 & -1 \\
                                       1 & 0 & 0 \\
                                       0 & 1 & 0 \\
                                     \end{pmatrix}, \ \ \ \lambda=\frac{3+\sqrt{5}}{2}=\mu^2.
$$
Let us denote $\frak{S}_V=\Phi(S_V)$. Then, away from $(0,0,0)$,
the family $\{\frak{S}_V\}$ is a smooth family of surfaces,
$\frak{S}_0$ is diffeomorphic to a cone, contains lines the
$\{y=0, z=0\}$ and $\{x=0, y=0\}$, and at each non-zero point on
these lines, it has a quadratic tangency with a horizontal or
vertical plane; compare Figures~\ref{fig:sing0.0} and
\ref{fig:sing0.2}.

We will use the variables $(x,y,z)$ for coordinates in
$\mathbb{R}^3$. For a point $p\in \mathbb{R}^3$, we will denote
its coordinates by $(x_p, y_p, z_p)$.

In order to study the properties of the map $f$ (i.e., of the map
$T$ in a small neighborhood of singularities), we need the
following statement.

\bprop\label{p.expanding} Given $C_1>0, C_2>0, \lambda>1$,
$\varepsilon \in (0,\frac{1}{4})$, and $\eta>0$, there exist
$\delta_0=\delta_0(C_1, C_2, \lambda, \varepsilon)$, $N_0\in
\mathbb{N}$, $N_0=N_0(C_1, C_2, \lambda, \varepsilon, \delta_0)$,
and $C=C(\eta)>0$ such that for
any $\delta\in (0, \delta_0)$, 
the following
holds.

Let $f:\mathbb{R}^3\to \mathbb{R}^3$ be a $C^2$-diffeomorphism such that
\begin{itemize}
\item[{\rm (i)}] $\|f\|_{C^2}\le C_1$; \item[{\rm (ii)}] The plane
$\{z=0\}$ is invariant under iterates of $f$; \item[{\rm (iii)}]
$\|Df(p)-A\|<\delta$ for every $p\in \mathbb{R}^3$, where
$$
A=\begin{pmatrix}
    \lambda^{-1} & 0 & 0 \\
    0 & 1 & 0 \\
    0 & 0 & \lambda \\
  \end{pmatrix}
$$
is a constant matrix.
\end{itemize}
Introduce the following cone field in $\mathbb{R}^3$:
\beq\label{e.coneformula} K_p=\{v\in T_p\mathbb{R}^3, \
\mathbf{v}=\mathbf{v}_{xy}+\mathbf{v}_z : |\mathbf{v}_z|\ge C_2
\sqrt{|z_p|}|\mathbf{v}_{xy}|\}. \eneq Given a  point $p=(x_p,
y_p, z_p)$ such that $0<z_p<1$, denote by $N=N(p)$ the smallest
integer $N\in \mathbb{N}$ such that $f^N(p)$ has $z$-coordinate
larger that 1.

If $N(p)\ge N_0$ {\rm (}i.e., if $z_p$ is small enough{\rm )} then
\beq\label{e.expansionforN} |Df^N(\mathbf{v})|\ge
\lambda^{\frac{N}{2}(1-4\varepsilon)}|\mathbf{v}| \ \ \ \text{for
any } \ \  \mathbf{v}\in K_p, \eneq and if
$Df^N(\mathbf{v})=\mathbf{u}=\mathbf{u}_{xy}+\mathbf{u}_z$, then
\beq\label{e.u} |\mathbf{u}_{xy}|<2\delta^{1/2}|\mathbf{u}_{z}|.
\eneq

Moreover, if $|\mathbf{v}_z|\ge \eta|\mathbf{v}_{xy}|$ then
\beq\label{e.expansionfork} |Df^k(\mathbf{v})|\ge
C\lambda^{\frac{k}{2}(1-4\varepsilon)}|\mathbf{v}| \ \ \ \text{for
each } \ \  k=1, 2, \ldots, N. \eneq \enprop

Returning to our specific situation at hand, if, for a given
$\delta$, the neighborhood $U$ is small enough, then at every
point $P\in U$ the differential $D(\Phi\circ T\circ \Phi^{-1})(P)$
satisfies the condition (iii) of Proposition \ref{p.expanding}.
Also, since the tangency of $\frak{S}_0$ with the horizontal plane
is quadratic, there exists $C_2>0$ such that every vector tangent
to $\frak{S}_0$ from the cone $D\Phi(\mathcal{K}^u)$ also belongs
to the cone \eqref{e.coneformula}. The same holds for vectors
tangent to $\frak{S}_V$ from continuations of cones
$D\Phi(\mathcal{K}^u)$ if $V$ is small enough. Therefore,
Proposition~\ref{p.expanding} can be applied to all those vectors.

We postpone the proof of Proposition~\ref{p.expanding} to
Section~\ref{s.ppe} and first show how to use it to prove Theorem
\ref{t.main}.

\section{Proof of Theorem~\ref{t.main} Assuming
Proposition~\ref{p.expanding}}\label{s.s05}

In order to prove the hyperbolicity of $\Omega_V$ we construct
only the unstable cone field and prove the unstable cone
condition. Due to the symmetry of the trace map, the stable cones
can be constructed in the same way.

Let $U'$ be a neighborhood of $\{P_1, P_2, P_3, P_4\}$ where the
results of the previous section can be applied. Since
$\mathbb{S}\backslash U'$ is compact, $F^{-1}(\mathbb{S}\backslash
U')$ is also compact. Denote
$$
\widetilde{C}=2\max_{p\in \mathbb{S}\backslash U'} \left\{
\|DF^{-1}(p)\|, \|DF(F^{-1}(p))\| \right\} < \infty.
$$
Take any $p\in \mathbb{S}\backslash U'$ and any $\mathbf{v}\in
T_p\mathbb{S}, \mathbf{v}\in \mathcal{K}^u.$ If $ T^n(p)\in
\mathbb{S}\backslash U'$, then \beq\label{e.ctilda}
\|DT^n(\mathbf{v})\|=\|D(F\circ \mathcal{A}^n\circ
F^{-1})(\mathbf{v})\|\ge \widetilde{C}^{-2} \mu^n \|\mathbf{v}\|.
\eneq Fix a small $\varepsilon>0$ and $n^*\in \mathbb{N}$ such
that $\widetilde{C}^{-2}\mu^{n^*}\ge \mu^{n^*(1-\varepsilon)}.$

\blm There exists a neighborhood $U^*\subset U'$ of the singular
set $\{P_1, P_2, P_3, P_4\}$ such that if $p\not \in U^*$ but
$T^{-1}(p)\in U^*$, and $n_0$ is the smallest positive integer
such that $T^{n_0}(p)\in U^*$, then the finite orbit $\{p, T(p),
\ldots, T^{n_0-1}(p)\}$ contains at least $n^*$ points outside of
$U'$. \elm
\begin{proof}
Take a point $p_1\in W^{uu}_{loc}(P_1)\subset U$, denote
$p_2=T(p_1)$, and consider the closed arc $J\subset W^{uu}(P_1)$
between points $p_1$ and $p_2$. For any point $p\in J$, denote by
$m(p)$ the smallest number $m\in \mathbb{N}$ such that the finite
orbit $\{p, T(p), \ldots, T^{m}(p)\}$ contains $n^*$ points
outside of $\overline{U'}$. Notice that the function $m(p)$ is
upper semi-continuous, and therefore (due to compactness of $J$)
bounded. Let $M\in \mathbb{N}$ be an upper bound. The set
$\bigcup_{i=0}^MT^i(J)$ is compact and does not contain the
singularities $P_1, P_2, P_3$ and $P_4$. Therefore, there exists
$\xi>0$ so small that if $\text{\rm dist}(q, J)<\xi$, then the
distance between any of the first $M$ iterates of $q$ and any of
the points $P_1, P_2, P_3, P_4$ is greater than $\xi$, and the
finite orbit $\{q, T(q), \ldots, T^{M}(q)\}$ contains at least
$n^*$ points outside of $\overline{U'}$. Now take $\xi'>0$ so
small that any point in $\xi'$-neighborhood of $P_1$ whose orbit
follows $W^{uu}(P_1)$ towards $p_1$ hits the $\xi$-neighborhood of
$J$ before leaving $U'$. Now we can take the $\xi'$-neighborhood
of the set $P_1, P_2, P_3, P_4$ as $U^*$.
\end{proof}

For small $V$, denote by $\mathbb{S}_{V, U^*}$ the bounded
component of $S_V\backslash U^*$. The family $\{\mathbb{S}_{V,
U^*}\}_{V\in [0,V_0)}$ of surfaces with boundary depends smoothly
on the parameter and has uniformly bounded curvature. For small
$V$, a projection $\pi_V:\mathbb{S}_{V, U^*} \to \mathbb{S}$ is
defined. The map $\pi_V$ is smooth, and if $p\in \mathbb{S}$,
$q\in \mathbb{S}_{V, U^*}$, and $\pi_V(q)=p$, then $T_p
\mathbb{S}$ and $T_q S_V$ are close. Denote by $\mathcal{K}_V^u$
(resp., $\mathcal{K}_V^s$) the image of the cone $\mathcal{K}^u$
(resp., $\mathcal{K}^s$) under the differential of $\pi_V^{-1}$.

Our choice of $n^*$ guarantees that the following statement holds.
 \blm\label{l.uprimeexpansion} There exists $V_0>0$ such that for
every $V\in [0, V_0)$, the following holds: if $\{q, T(q), T^2(q),
\ldots, T^n(q)\}\subset \mathbb{S}_{V, U^*}$, $\mathbf{v}\in
\mathcal{K}_V^u\subset T_qS_V$, $q, T^n(q)\in \mathbb{S}_{V, U'}$,
and $n\ge n^*$, then
 $$
 \|DT^n(\mathbf{v})\|\ge \mu^{n(1-2\varepsilon)}\|\mathbf{v}\|.
 $$
 \elm
 \blm\label{l.exansionlong}
There exist $V_0>0$ and $C>0$ such that for any $V\in [0, V_0)$,
we have that if $q\in \mathbb{S}_{V, U'}$, $\mathbf{v}\in
\mathcal{K}_V^u\subset T_qS_V$, and $T^n(q)\in \mathbb{S}_{V,
U'}$, then
 $$
 \|DT^n(\mathbf{v})\|\ge C\mu^{n(1-4\varepsilon)}\|\mathbf{v}\|.
 $$
 \elm
 \begin{proof}
 Let us split the orbit $\{q, T(q), T^2(q), \ldots, T^n(q)\}$ into several intervals
 $$
 \{q, T(q), T^2(q), \ldots, T^{k_1-1}(q)\}, \{T^{k_1}(q),  \ldots, T^{k_2-1}(q)\}, \ldots, \{T^{k_s}(q), \ldots, T^n(q)\}
 $$
 in such a way that the following properties hold:
 \begin{enumerate}
 \item
 for each $  i=1, 2, \ldots, s$, the points $T^{k_i-1}(q)$ and $T^{k_i}(q)$ are outside of $U'$;
 \item
 if $\{T^{k_i}(q),  \ldots, T^{k_{i+1}-1}(q)\}\cap U^*\ne \emptyset$, then $\{T^{k_i+1}(q),  \ldots, T^{k_{i+1}-2}(q)\}\subset U'$;
 \item
 for each $  i=1, 2, \ldots, s-1$, we have either $k_{i+1}-k_i\ge n^*$ or $\{T^{k_i}(q),  \ldots, T^{k_{i+1}-1}(q)\}\cap U^*\ne \emptyset$.
 \end{enumerate}
 Such a splitting exists due to the choice of $U^*$ above.

The following lemma is a consequence of Lemma \ref{l.differential}
and property~\eqref{e.u} from Proposition \ref{p.expanding}.
 \blm\label{l.coneinvariance}
Suppose $q\in \mathbb{S}_{V, U'}$, $\mathbf{v}\in
\mathcal{K}_V^u\subset T_qS_V$, $T(q)\in U'$, and $l\in
\mathbb{N}$ is the smallest number such that $T^l(q)\not \in U'$.
If $l$ is large enough and $V$ is small enough, then
 $$
 DT^l(q)(\mathbf{v})\in \mathcal{K}_V^u\subset T_{T^l(q)}S_V.
 $$
 \elm
Apply Lemma \ref{l.coneinvariance} to those intervals in the
splitting that are contained in $U'$.  Notice that largeness of
$l$ can be provided by the choice of $U^*$. Together with
Proposition \ref{p.expanding} applied to these intervals, and
Lemma \ref{l.uprimeexpansion} applied to intervals that do not
intersect $U^*$, this guarantees uniform expansion of
$\mathbf{v}$. The first and the last interval may have length
greater than $n^*$, and then Lemma \ref{l.uprimeexpansion} can be
applied, or smaller than $n^*$, but then taking a small enough
constant $C$ (say, $C<\widetilde{C}^{-2}$) will compensate for the
lack of uniform expansion on these intervals. Lemma
\ref{l.exansionlong} is proved.
 \end{proof}

For every small $V>0$, there exists $\eta_V>0$ ($\eta_V\to 0$ as
$V\to 0$) such that if $p\not \in U'$, $T(p)\in U'$, and
$\mathbf{v}\in \mathcal{K}_V^u(p)$, then for the vector
$\mathbf{w}\equiv D\Phi^{-1}(p)(\mathbf{v})$,
$\mathbf{w}=\mathbf{w}_z+ \mathbf{w}_{xy}$, we have
$|\mathbf{w}_z|\ge \eta_V|\mathbf{w}_{xy}|$. An application of
Proposition \ref{p.expanding} together with Lemma
\ref{l.exansionlong} proves uniform expansion of vectors from
$\mathcal{K}_V^u$. Due to the symmetries of the trace map, all
vectors from $\mathcal{K}_V^s$ are expanded by iterates of
$T^{-1}$. Thus, hyperbolicity of the set $\Omega_V$ follows and
Theorem~\ref{t.main}~(i) is proved.

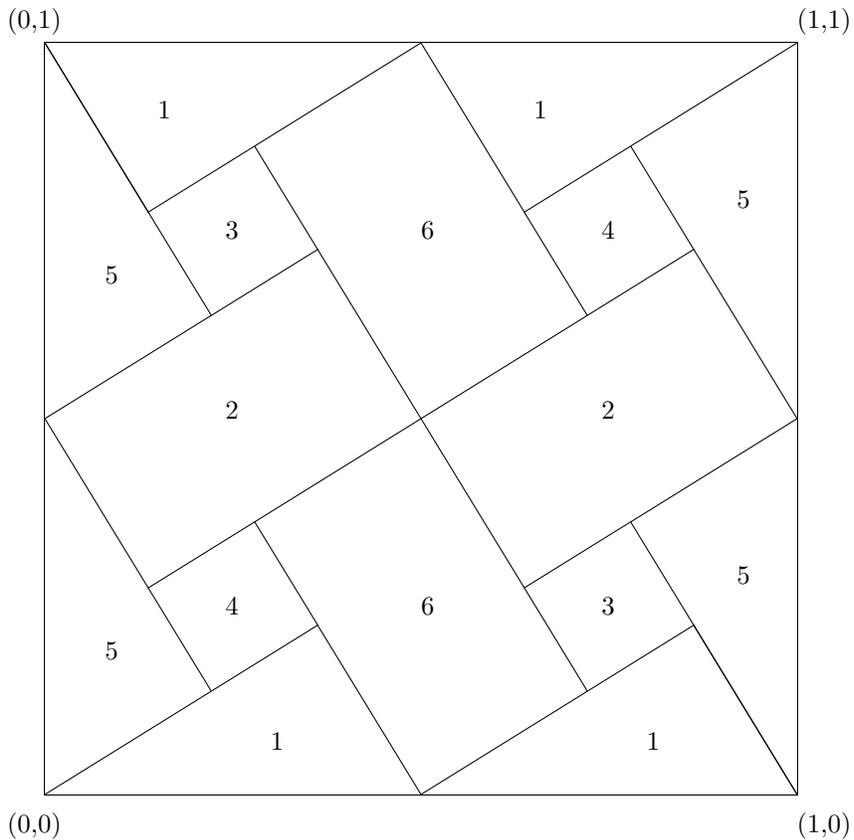
\begin{figure}
\setlength{\unitlength}{1mm}
\begin{picture}(120,120)

\put(10,10){\framebox(100,100)}

\put(5,5){(0,0)}

\put(110,5){(1,0)}

\put(5,112){(0,1)}

\put(110,112){(1,1)}

\put(10,10){\line(161,100){36.2}}

\put(60,10){\line(-61,100){13.8}}

\put(40,16){$1$}

\put(60,10){\line(161,100){36.2}}

\put(110,10){\line(-61,100){13.8}}

\put(90,16){$1$}

\put(60,110){\line(-161,-100){36.2}}

\put(10,110){\line(61,-100){13.8}}

\put(25,100){$1$}

\put(110,110){\line(-161,-100){36.2}}

\put(60,110){\line(61,-100){13.8}}

\put(75,100){$1$}

\put(10,60){\line(61,-100){22.1}}

\put(18,28){$5$}

\put(10,60){\line(161,100){22.1}}

\put(10,110){\line(61,-100){22.1}}

\put(18,78){$5$}

\put(110,60){\line(-61,100){22.1}}

\put(102,38){$5$}

\put(110,60){\line(-161,-100){22.1}}

\put(110,10){\line(-61,100){22.1}}

\put(102,88){$5$}

\put(60,60){\line(-161,-100){36.3}}

\put(60,60){\line(161,100){36.3}}

\put(60,60){\line(-61,100){22.1}}

\put(60,60){\line(61,-100){22.1}}

\put(88,46.35){\line(-161,-100){14.3}}

\put(32,73.65){\line(161,100){14.3}}

\put(37.88,46.28){\line(61,-100){8.4}}

\put(82.12,73.72){\line(-61,100){8.4}}

\put(34,34){$4$}

\put(84,84){$4$}

\put(34,84){$3$}

\put(84,34){$3$}

\put(34,60){$2$}

\put(84,60){$2$}

\put(60,34){$6$}

\put(60,84){$6$}

\end{picture}
\caption{The Markov partition for
$T|_{\mathbb{S}}$.}\label{fig.Casdagli-Markov}
\end{figure}

The Markov partition for $T|_{\mathbb{S}}$ (i.e., for $V=0$) was
presented explicitly by Casdagli in~\cite[Section~2]{Cas}; compare
Figure~\ref{fig.Casdagli-Markov}.

Since the Markov partition is formed by finite pieces of  strong
stable and strong unstable manifolds of the periodic points of $T$
and these manifolds depend smoothly on the parameter $V$, there
exists a Markov partition for $\Omega_V$ with the same matrix as
for $V=0$ (see \cite[Appendix 2]{PT} for more details on Markov
partitions for 2-dimensional hyperbolic maps). This establishes
Theorem~\ref{t.main}~(ii).

\blm\label{l.transversality} If $V$ is small enough, then the line
$\ell_V$ is transversal to the cone field $\mathcal{K}_V^s$. \elm

\begin{proof}
For $V=0$ and small enough $\zeta$ in \eqref{e.cones}, this is
true since $F^{-1}(l_0)=\{\theta=-\varphi\}$, and the vector $(1,
-1)$ is not an eigenvector of $A$. Therefore this is also true for
every sufficiently small $V$ by continuity.
\end{proof}

In order to show that $\ell_V$ is also transversal to the stable
manifolds of $\Omega_V$ inside of $U^*$, let us consider the rectifying
coordinates $\Phi:U^*\to \mathbb{R}^3$ again and define the
central-unstable cone field in $\Phi(U^*)$:
$$
K_p^{cu}=\left\{\mathbf{v}\in T_p\mathbb{R}^3,
\mathbf{v}=\mathbf{v}_x+\mathbf{v}_{yz} :
|\mathbf{v}_{yz}|>\zeta^{-1}|\mathbf{v}_x|\right\}.
$$
Since $\Phi(l_0)$ is transversal to the plane $\{x=0\}$, the curve
$\Phi(\ell_V)$ is transversal to this invariant cone field if
$\zeta, V$, and $U^*$ are small enough. Every stable manifold of
$\Omega_V$ in the rectifying coordinates is tangent to this
central-unstable cone field, and together with Lemma
\ref{l.transversality} this implies that $\ell_V$ is transversal
to stable manifolds of $\Omega_V$. This shows
Theorem~\ref{t.main}~(iii) and hence concludes the proof of
Theorem~\ref{t.main}.

\section{Proof of Proposition~\ref{p.expanding}}\label{s.ppe}

\subsection{Properties of the Recurrent Sequences}

In this subsection we formulate and prove several lemmas on
recurrent sequences that will be used in the next subsection to
prove Proposition~\ref{p.expanding}.

\blm\label{l.dD} Given $C_1>0$, $C_2>0$, $\lambda>1$, and
$\varepsilon \in (0,\frac{1}{4})$, there exist
$\delta_0=\delta_0(C_1, C_2, \lambda, \varepsilon)$ and $N_0\in
\mathbb{N}$, $N_0=N_0(C_1, C_2, \lambda, \varepsilon, \delta_0)$,
such that for every $\delta\in (0, \delta_0)$ and every $N\ge
N_0$, the following holds. Suppose that the sequences
$\{d_i\}_{i=0}^N$, $\{D_i\}_{i=0}^N$ are defined by the initial
conditions \beq\label{e.initialDd} \left\{
  \begin{array}{ll}
    d_0=1 \\
    D_0\ge C_2\frac{1}{(\lambda+\delta)^{N/2}},
  \end{array}
\right.
\eneq
and recurrence relations
\beq\label{e.recurrenceDd}
\left\{
  \begin{array}{ll}
    d_{k+1}=(1+2\delta)d_k+\delta D_k, \\
    D_{k+1}=(\lambda-\delta)D_k-C_1b_kd_k,
  \end{array}
\right.
\eneq
where
$$
b_k=(\lambda-\delta)^{-N+k}.
$$
Then
$$
\left\{
  \begin{array}{ll}
    d_N\le 2\delta^{1/2}D_N \\
    D_N\ge
D_0\lambda^{N(1-\varepsilon)}>\lambda^{\frac{N}{2}(1-4\varepsilon)}.
  \end{array}
\right.
$$
 \elm
\brm\label{r.coneDd} Notice that Lemma \ref{l.dD} implies also that $$
\frac{d_N}{D_N}\le 2\delta^{1/2}.
$$
This inequality will allow us to obtain a small cone (of size of
order $\delta^{1/2}$) where the image of a vector after leaving a
neighborhood of a singularity is located. \erm

\begin{proof}[Proof of Lemma~\ref{l.dD}.]
Let us denote
$$
\Lambda_{-}=\lambda^{1-\varepsilon},\ \ \
\Lambda_{+}=\lambda^{1+\varepsilon}.
$$
If $\delta$ is small enough, then
$1<\Lambda_{-}<\lambda-\delta<\lambda+\delta<\Lambda_{+}$. We will
prove by induction that \beq\label{e.inductiveassumptions} \left\{
  \begin{array}{ll}
    d_k\le (1+2\delta+\delta^{1/2})\max\{d_{k-1}, \delta^{1/2} D_{k-1}\} \\
    D_k\ge \Lambda_{-}D_{k-1}\ge  \Lambda_{-}^k D_0.
  \end{array}
\right. \eneq
It is clear that these inequalities for $k=N$ imply Lemma \ref{l.dD}.

Let us first check the base of induction. If $N$ is large, then
$D_0<1$, so $\max\{d_0, \delta^{1/2}D_0\}=d_0=1$. So we have
 $$
 d_1=(1+2\delta)d_0+\delta D_0<1+3\delta<(1+2\delta+\delta^{1/2})d_0.
 $$
 This is the first inequality in (\ref{e.inductiveassumptions}) for $k=1$.

 Also we have
\begin{align*}
D_1 & = (\lambda-\delta)D_0-C_1b_0d_0 \\
& = (\lambda-\delta)D_0\left(1-\frac{C_1b_0d_0}{D_0}\right) \\
& \ge (\lambda-\delta) D_0 \left(1-\frac{C_1(\lambda+\delta)^{N/2}}{C_2(\lambda-\delta)^{N}}\right)\\
& > (\lambda - \delta) D_0 \left( 1 - \frac{C_1\Lambda_{+}^{N/2}}{C_2\Lambda_{-}^N} \right) \\
& =( \lambda - \delta ) D_0 \left( 1 - \frac{C_1}{C_2}
\lambda^{-\frac{N}{2} ( 1 - 3 \varepsilon)} \right),
\end{align*}
and since
$(\lambda-\delta)\left(1-\frac{C_1}{C_2}\lambda^{-\frac{N}{2}(1-3\varepsilon)}\right)>\Lambda_{-}$
if $N$ is large enough, we have $D_1>\Lambda_{-}D_0$. We checked
the base of induction.

We proceed to proving the induction step. Assume that for some
$k$, the inequalities (\ref{e.inductiveassumptions}) hold. Let us
show that these inequalities also hold for $k+1$. We have
$$
d_{k+1}=(1+2\delta)d_k+\delta D_k.
$$
If $\delta^{1/2}D_k\le d_k$, then
\begin{align*}
d_{k+1} & \le ( 1 + 2 \delta ) d_k + \delta^{1/2} d_k \\
& = ( 1 + 2 \delta + \delta^{1/2} ) d_k \\
& = ( 1 + 2 \delta + \delta^{1/2}) \max\{d_{k-1}, \delta^{1/2}
D_{k-1}\}.
\end{align*}
If $\delta^{1/2}D_k> d_k$, then
\begin{align*}
d_{k+1} & \le (1+2\delta)\delta^{1/2}D_k+\delta D_k \\
& = \delta^{1/2} D_k ( 1 + 2 \delta + \delta^{1/2} ) \\
& = ( 1 + 2 \delta + \delta^{1/2})\max\{d_{k-1}, \delta^{1/2}
D_{k-1}\}.
\end{align*}

In order to estimate $D_{k+1}$, we need more information on
$\{d_i\}_{i=0}^k$. \blm\label{l.dsubk} Assume that $k^*\le N$ and
\eqref{e.inductiveassumptions} holds for all $k=1, \ldots, k^*$.
Assume also that for $k=1, 2, \ldots, l-1$, we have
$\delta^{1/2}D_k\le d_k$ and $\delta^{1/2}D_l>d_l$. Then,
$\delta^{1/2}D_k>d_k$ for every $k=l, l+1, \ldots, k^*$. \elm
\begin{proof}
If $\delta$ is small enough, $1+2\delta+\delta^{1/2}<\Lambda_{-}$.
Therefore, if $\delta^{1/2}D_l>d_l$, then
\begin{align*}
d_{l+1} & = (1+2\delta)d_l+\delta D_l \\
& < (1 + 2\delta+\delta^{1/2})\delta^{1/2} D_l \\
& < \Lambda_{-}\delta^{1/2}D_l \\
& \le \delta^{1/2} D_{l+1}.
\end{align*}
In the same way, $d_{l+2}<\delta^{1/2}D_{l+2}$, and so on.
\end{proof}
Notice that Lemma \ref{l.dsubk} immediately implies the following
statement. \blm\label{l.dsubk1} Assume that $k^*\le N$ and
\eqref{e.inductiveassumptions} holds for all $k=1, \ldots, k^*$.
If $d_k\ge \delta^{1/2}D_k$ for some $k\in\{1, 2, \ldots, k^*\}$,
then $d_k\le (1+2\delta +\delta^{1/2})^k$. \elm

Now let us estimate $D_{k+1}$. If $d_k<\delta^{1/2}D_k$, then
\begin{align*}
D_{k+1} & = (\lambda-\delta)D_k
\left(1-\frac{C_1b_kd_k}{(\lambda-\delta)D_k} \right) \\
& \ge (\lambda-\delta)D_k \left(1 -
\left(\frac{C_1}{\lambda-\delta} \right)
\frac{\delta^{1/2}}{(\lambda-\delta)^{N-k}}\right) \\
& > (\lambda-\delta)D_k(1-C_1\delta^{1/2}) \\
& > \Lambda_{-}D_k
\end{align*}
if $\delta$ is small enough.

If $d_k \ge \delta^{1/2}D_k$, then due to Lemma \ref{l.dsubk1}, we
have $d_k\le (1+2\delta+\delta^{1/2})^k$, and hence
\begin{align*}
D_{k+1} & = (\lambda-\delta)D_k\left(1-\frac{C_1b_kd_k}{(\lambda-\delta)D_k}\right) \\
& \ge ( \lambda - \delta ) D_k \left( 1 - \left(\frac{C_1}{\lambda-\delta}\right)
\frac{(1 + 2 \delta + \delta^{1/2})^k}{\Lambda_{-}^k D_0(\lambda-\delta)^{N-k}}\right)\\
& \ge ( \lambda - \delta ) D_k \left(1 - \left(
\frac{C_1}{C_2(\lambda-\delta)} \right) \frac{(1 + 2 \delta + \delta^{1/2})^k
(\lambda + \delta)^{N/2}}{\Lambda_{-}^k(\lambda - \delta)^{N-k}}\right)\\
& \ge ( \lambda - \delta) D_k \left(1 - \left(\frac{C_1}{C_2(\lambda-\delta)}\right)
\frac{(1 + 2\delta+\delta^{1/2})^k\Lambda_{+}^{N/2}}{\Lambda_{-}^k\Lambda_{-}^{N-k}}\right)\\
& \ge ( \lambda - \delta) D_k \left(1 - \left(\frac{C_1}{C_2(\lambda-\delta)}\right)
\frac{(1+2\delta+\delta^{1/2})^k \lambda^{\frac{N}{2}(1+\varepsilon)}}{\lambda^{N(1-\varepsilon)}}\right) \\
& \ge ( \lambda - \delta ) D_k \left(1 -
\left(\frac{C_1}{C_2(\lambda-\delta)}\right)
(1+2\delta+\delta^{1/2})^k
\lambda^{-\frac{N}{2}(1-3\varepsilon)}\right)\\
& >\Lambda_{-}D_k
\end{align*}
if $N$ is large enough. This concludes the proof of Lemma
\ref{l.dD}.
\end{proof}

\blm\label{l.Aa} Given $C_1>0$, $C_2>0$, $\lambda>1$, and
$\varepsilon \in (0,\frac{1}{4})$, there exist
$\delta_0=\delta_0(C_1, C_2, \lambda, \varepsilon)$ and $N_0\in
\mathbb{N}$, $N_0=N_0(C_1, C_2, \lambda, \varepsilon, \delta_0)$,
such that for every $\delta\in (0, \delta_0)$ and every $N\ge
N_0$, the following holds. Assume that the sequences
$\{a_k\}_{k=0}^N$ and $\{A_k\}_{k=0}^N$ have the following
properties: \beq\label{e.initialAa} \left\{
  \begin{array}{ll}
    a_0=1 \\
    A_0\ge C_2\sqrt{\widetilde{b}_0}>0,
  \end{array}
\right.
\eneq
and for $k=0, 2, \ldots, N-1$
\beq\label{e.recurrenceAa}
\left\{
  \begin{array}{ll}
    a_{k+1}\le (1+2\delta)a_k+\delta A_k, \\
    A_{k+1}\ge (\lambda-\delta)A_k-C_1\widetilde{b}_k a_k,
  \end{array}
\right.
\eneq
where the sequence $\widetilde{b_k}$ has the properties:
$$
0<\widetilde{b}_0<\widetilde{b}_1 <\ldots <
\widetilde{b}_{N-1}<1\le \widetilde{b}_N ,
$$
and for all $k=0, 2, \ldots, N-1$
$$
(\lambda-\delta)\widetilde{b}_k \le \widetilde{b}_{k+1} \le
(\lambda+\delta)\widetilde{b}_k.
$$
Then,
$$
\left\{
  \begin{array}{ll}
    a_N\le 2\delta^{1/2}A_N \\
    A_N\ge \lambda^{N(1-\varepsilon)}A_0> \lambda^{\frac{N}{2}(1-4\varepsilon)}.
  \end{array}
\right.
$$
\elm
\begin{proof}
Consider the sequences $\{d_k\}$ and $\{D_k\}$ defined by
(\ref{e.initialDd}) and (\ref{e.recurrenceDd}), where we take
$D_0=A_0$. If we prove that for all $k=0, 1, \ldots, N$
\beq\label{e.AaDdinduction} \left\{
  \begin{array}{ll}
    A_k\ge D_k \\
    \frac{A_k}{a_k}\ge \frac{D_k}{d_k},
  \end{array}
\right.
\eneq
then the required inequalities follow from Lemma \ref{l.dD}.

We will prove (\ref{e.AaDdinduction}) by induction. The base of
induction ($k=0$) is provided by our choice of $D_0$. Let us turn
to the induction step. If \eqref{e.AaDdinduction} holds for some
$k$, then
\begin{align*}
\frac{A_{k+1}}{a_{k+1}} & \ge
\frac{(\lambda-\delta)A_k-C_1\widetilde{b}_k
a_k}{(1+2\delta)a_k+\delta
A_k} \\
& =(\lambda-\delta)\delta^{-1}-\frac{(1+2\delta)(\lambda-\delta)\delta^{-1}+C_1\widetilde{b}_k}{(1+2\delta)+\delta\frac{A_k}{a_k}}\\
& \ge (\lambda - \delta) \delta^{-1} - \frac{(1 + 2\delta)(\lambda
- \delta)\delta^{-1}+C_1{b_k}}{(1+2\delta)+\delta\frac{D_k}{d_k}} \\
& =  \frac{D_{k+1}}{d_{k+1}}.
\end{align*}
Also,
\begin{align*}
A_{k+1} & \ge (\lambda-\delta)A_k-C_1b_ka_k \\
& = A_k\left((\lambda-\delta)-C_1\widetilde{b}_k \frac{a_k}{A_k}\right) \\
& \ge D_k\left((\lambda-\delta)-C_1{b_k}\frac{d_k}{D_k}\right) \\
& = D_{k+1}.
\end{align*}
 Lemma \ref{l.Aa} is proved.
\end{proof}

\subsection{Expansion of Vectors From Large Cones}

Finally, we come to the

\begin{proof}[Proof of Proposition~\ref{p.expanding}.]
Take $\mathbf{v}\in K_p$. If $|\mathbf{v}_z|>|\mathbf{v}_{xy}|$,
then the required inequalities follow just from the condition
(iii). Therefore, we can assume that $|\mathbf{v}_{xy}|\ge
|\mathbf{v}_z|$, and normalize the vector $\mathbf{v}$ assuming
that $|\mathbf{v}_{xy}|=1$.

Denote
$$
Df(p)=\begin{pmatrix}
        \nu(p) & m_1(p) & t_1(p) \\
        m_2(p) & e(p) & t_2(p) \\
        s_1(p) & s_2(p) & \lambda(p) \\
      \end{pmatrix}.
$$
Given a vector $\mathbf{v}\in K_p$, denote
$\mathbf{v}'=Df(p)(\mathbf{v})$,
$\mathbf{v}'=\mathbf{v}'_{xy}+\mathbf{v}'_z$,
$\mathbf{v}'_{xy}=\mathbf{v}'_x+\mathbf{v}'_y$. We have
$$
Df(p)(\mathbf{v})=\begin{pmatrix}
        \nu(p) & m_1(p) & t_1(p) \\
        m_2(p) & e(p) & t_2(p) \\
        s_1(p) & s_2(p) & \lambda(p) \\
      \end{pmatrix}\begin{pmatrix}
                     \mathbf{v}_x \\
                     \mathbf{v}_y \\
                     \mathbf{v}_z \\
                   \end{pmatrix}=\begin{pmatrix}
                     \nu(p)\mathbf{v}_x+m_1(p)\mathbf{v}_y+t_1(p)\mathbf{v}_z \\
                     m_2(p)\mathbf{v}_x+e(p)\mathbf{v}_y+t_2(p)\mathbf{v}_z \\
                     s_1(p)\mathbf{v}_x+s_2(p)\mathbf{v}_y+\lambda(p)\mathbf{v}_z \\
                   \end{pmatrix}
$$
The conditions (i)--(iii) of Proposition \ref{p.expanding} imply
that $|\nu(p)|\le \lambda^{-1}+\delta$, $|m_1(p)|, |m_2(p)|,
|t_1(p)|, |t_2(p)|\le \delta$, and $|\lambda(p)|\ge
\lambda-\delta$. Also, if $p$ belongs to the plane $\{z=0\}$, then
$s_1(p)=s_2(p)=0$. Since $\|f\|_{C^2}\le C_1$, for arbitrary $p$,
we have $|s_1(p)|, |s_2(p)|\le C_1z_p$.

Let us use the norm $\|\mathbf{v}\|=|\mathbf{v}_x|+|\mathbf{v}_y|+|\mathbf{v}_z|$.
Then we have
$$
\|\mathbf{v}_{xy}'\|\le
(1+2\delta)\|\mathbf{v}_{xy}\|+\delta\|\mathbf{v}_z\|,
$$
$$
\|\mathbf{v}_z'\|\ge
(\lambda-\delta)\|\mathbf{v}_z\|-C_1z_p\|\mathbf{v}_{xy}\|.
$$
Denote by $\widetilde{b}_k$ the $z$-coordinate of $f^k(p)$, by
$A_k$ the $z$-component of the vector $Df^k(p)(\mathbf{v})$, and
by $a_k$ the $xy$-component of $Df^k(p)(\mathbf{v})$. Applying
Lemma \ref{l.Aa} we get \eqref{e.expansionforN}, and \eqref{e.u}
follows from Remark~\ref{r.coneDd}.

Now if $A_0=|\mathbf{v}_z|\ge \eta$, then (applying
Lemma~\ref{l.dD} for $D_0=\eta$) from \eqref{e.AaDdinduction} and
\eqref{e.inductiveassumptions}, we have
$$
|Df^k(p)(\mathbf{v})|\ge A_k\ge D_k\ge \Lambda_{-}^kD_0\ge
\lambda^{\frac{k}{2}(1-\varepsilon)}\eta>\frac{1}{2}\lambda^{\frac{k}{2}(1-\varepsilon)}\eta
|\mathbf{v}|>C\lambda^{\frac{k}{2}(1-4\varepsilon)} |\mathbf{v}|
$$
for every $k=1, 2, \ldots, N$, where $C=\frac{1}{2}\eta$. Proposition \ref{p.expanding} is proved.
\end{proof}

\end{document}